\documentclass[11pt]{amsart}   

\usepackage{graphics} 
\usepackage{amscd}    

\newtheorem{theorem}{Theorem}[section]
\newtheorem{lemma}[theorem]{Lemma}
\newtheorem{proposition}[theorem]{Proposition}
\newtheorem{corollary}[theorem]{Corollary}
\newtheorem{question}[theorem]{Question}
\theoremstyle{definition}

\theoremstyle{remark}
\newtheorem{remark}[theorem]{Remark}
\theoremstyle{definition}
\newtheorem{example}[theorem]{Example}

\numberwithin{equation}{section}

\usepackage{amssymb}
\usepackage[leqno]{amsmath}
\usepackage{amsfonts}
\usepackage{amsopn}
\usepackage{amstext}
\usepackage{amsthm}

\providecommand{\cal}{\mathcal}
\renewcommand{\Bbb}{\mathbb}

\newenvironment{pf}{\begin{proof}}{\end{proof}}

\newcommand{\Aaa}{{\cal{A}}}

\newcommand{\Yu}{{\cal{U}}}
\newcommand{\Vee}{{\cal{V}}}

\newcommand{\Qyu}{{\Bbb{Q}}}
\newcommand{\Err}{{\Bbb{R}}}

\newcommand{\al}{\alpha}

\newcommand{\sig}{\sigma}
\newcommand{\eps}{\varepsilon}
\renewcommand{\phi}{\varphi}
\renewcommand{\rho}{\varrho}

\newcommand{\rest}{\restriction}

\newcommand{\ntr}{n\in\omega}

\newcommand{\loe}{\leqslant}
\newcommand{\goe}{\geqslant}

\newcommand{\subs}{\subseteq}

\newcommand{\nnempty}{\ne\emptyset}

\renewcommand{\iff}{\Longleftrightarrow}

\newcommand{\cl}{\operatorname{cl}}
\newcommand{\Int}{\operatorname{int}}
\newcommand{\w}{\operatorname{w}}

\newcommand{\bbl}{{\mathbb L}}
\newcommand{\bbk}{{\mathbb K}}
\newcommand{\ult}{\operatorname{Ult}}

\newcommand{\id}{\operatorname{id}}

\newcommand{\supp}{\operatorname{suppt}}

\newcommand{\pr}{\operatorname{pr}}
\newcommand{\liminv}{\varprojlim}

\newcommand{\poset}{{\Bbb{P}}}

\newcommand{\Es}{{\mathbb{S}}}

\newcommand{\Land}{\;\&\;}

\newcommand{\concat}{{}^\smallfrown}

\newcommand{\setof}[2]{\{#1\colon #2\}}

\newcommand{\sett}[2]{\{#1\}_{#2}}
\newcommand{\sn}[1]{\{#1\}} 
\newcommand{\dn}[2]{\{#1,#2\}} 
\newcommand{\pair}[2]{\langle #1, #2 \rangle} 
\newcommand{\map}[3]{#1\colon #2 \to #3} 
\newcommand{\img}[2]{#1[#2]} 
\newcommand{\inv}[2]{{#1}^{-1}[#2]} 

\newcommand{\suppt}{\supp}

\newcommand{\Ssig}{{\sig\in\Sigma}}

\newcommand{\invsys}[5]{\langle {#1}_{#4};{#2}_{#4}^{#5};#3 \rangle}
\renewcommand{\S}{\mathbb S}
\newcommand{\diag}{\Delta} 

\newcommand{\En}{{\mathcal N}}

\begin{document}
\noindent                                             
\begin{picture}(150,36)                               
\put(5,20){\tiny{Submitted to}}                       
\put(5,7){\textbf{Topology Proceedings}}              
\put(0,0){\framebox(140,34){}}                        
\put(2,2){\framebox(136,30){}}                        
\end{picture}                                         

\vspace{0.5in}

\title{Semi-Eberlein spaces}

\author{Wies{\l}aw Kubi\'s}
\address{Institute of Mathematics, University of Silesia, ul. Bankowa 14, 40-007 Katowice, Poland}
\curraddr{Department of Mathematics, York University,
Toronto, Canada}
\email{kubis@ux2.math.us.edu.pl, wkubis@yorku.ca}

\author{Arkady Leiderman}
\address{Department of Mathematics, Ben Gurion University of the Negev,
Be'er Sheva, Israel}
\email{arkady@math.bgu.ac.il}

\subjclass{Primary: 54D30; Secondary: 54C35, 54C10.}

\keywords{Semi-Eberlein space, Valdivia/Corson compact, semi-open retraction}


\begin{abstract} We investigate the class of compact spaces which are embeddable into a power of the real line $\Err^\kappa$ in such a way that $c_0(\kappa)=\setof{f\in\Err^\kappa}{(\forall\;\eps>0)\;|\setof{\al\in\kappa}{|f(\al)|>\eps}|<\aleph_0}$ is dense in the image. We show that this is a proper subclass of the class of Valdivia, even when restricted to Corson compacta. We prove a preservation result concerning inverse sequences with semi-open retractions. As a corollary we obtain that retracts of Cantor or Tikhonov cubes belong to the above class.
\end{abstract}

\maketitle

\section{Introduction}

This study is motivated by results on the class of Valdivia compact spaces, i.e. compact spaces embeddable into $\Err^\kappa$ in such a way that the $\Sigma$-product $\Sigma(\kappa)=\setof{f\in\Err^\kappa}{|\suppt(f)|\loe\aleph_0}$ is dense in the image. This class was first considered by Argyros, Mercourakis and Negrepontis in \cite{AMN} (the name {\it Valdivia compact} was introduced in \cite{DG}) and then studied by several authors, see e.g. \cite{Valdivia1990, Valdivia1991, DG, Kalenda4, Kalenda5, KU, KM}. Several important results on Valdivia compacta were established by Kalenda and we refer to his survey article \cite{Kalenda} for further references. Probably the most remarkable property of a Valdivia compact space is the existence of ``many retractions", namely if $K$ is Valdivia compact and $\map fKY$ is a continuous map then there exists a ``canonical" retraction $\map rKK$ such that $\img rK$ is  Valdivia of the same weight as $\img fK$ and $f=fr$ holds. This property implies that the Banach space $C(K)$ has a {\it projectional resolution of the identity\/}, a useful property from which one can deduce e.g. the existence of an equivalent locally uniformly convex norm (see \cite[Chapter VII]{DGZ} or \cite[Chapter 6]{F}). 

In this note we consider compact spaces which have a better embedding than Valdivia: namely, compact spaces embeddable into $\Err^\kappa$ in such a way that $c_0(\kappa)$ is dense in the image. We call such spaces {\em semi-Eberlein}. This class of spaces clearly contains all Eberlein compacta. We show that not all Corson compacta are semi-Eberlein.

An obstacle in proving results parallel to the class of Valdivia compacta is the fact that, in contrast to $\Sigma$-products, $c_0(\kappa)$ is not countably compact. As a rule, obvious modifications of arguments concerning Valdivia compacta cannot be applied to semi-Eberlein spaces.

We state some simple characterizations of semi-Eberlein compacta involving families of open $F_\sig$ sets and properties of spaces of continuous functions with the topology of pointwise convergence on a dense set. We prove a preservation theorem involving inverse sequences whose bonding maps are semi-open retractions with a metrizable kernel.
We show that semi-Eberlein compacta do not contain P-points. This immediately implies that the linearly ordered space $\omega_1+1$ is not semi-Eberlein. Finally, we show that the Corson compact space of Todor\v{c}evi\'c \cite{To2} (see also \cite{To1}) 
is not semi-Eberlein.
We finish by formulating several open questions.

\subsection*{Acknowledgements}
We would like to thank Stevo Todor\v{c}evi\'c for providing the argument in Example \ref{bajer}. We are grateful to the anonymous referee for his valuable remarks which helped us improve the presentation.

\section{Notation and definitions}
All topological spaces are assumed to be completely regular. By a ``map" we mean a continuous map.
Fix a set $S$. As mentioned in the introduction, $\Sigma(S)$ denotes the $\Sigma$-product of $S$ copies of $\Err$, i.e. the set of all functions $x\in\Err^S$ such that $\suppt(x):=\setof{s\in S}{x(s)\ne0}$ (the support of $x$) is countable. 
$c_0(S)$ denotes the subspace of $\Sigma(S)$ consisting of all functions $x\in\Err^S$ such that for every $\eps>0$ the set $\setof{s\in S}{|x(s)|>\eps}$ is finite.
Given $T\subs S$ we denote by $\pr_T$ the canonical projection from $\Err^S$ onto $\Err^T$. Note that $\img{\pr_T}{\Sigma(S)}=\Sigma(T)$ and $\img{\pr_T}{c_0(S)}=c_0(T)$.

Given two maps $\map fXY$ and $\map gXZ$ we denote by $f\diag g$ their diagonal product, i.e. $\map {f\diag g}X{Y\times Z}$ is defined by $(f\diag g)(x) = \pair{f(x)}{g(x)}$. 

A compact space $K$ is {\em Eberlein compact} if $K\subs c_0(S)$ for some $S$ or equivalently $K$ is embeddable into a Banach space with the weak topology (by Amir, Lindenstrauss \cite{AL}). 

Let $T$ be a set, a family $\Aaa$ of subsets $T$ is called {\em adequate} \cite{Ta1} provided
1) $\{t\}\in\Aaa$ for each $t\in T$; and 2) $A\in\Aaa$ iff $M\in \Aaa$ for any finite $M\subs A$. Every adequate family is a compact 0-dimensional space, when identifying sets with their characteristic functions lying in the Cantor cube $\{0,1\}^T$. Any space homeomorphic to an adequate family of sets is called {\em adequate compact}.
A collection of sets $\Yu$ is {\em $T_0$ separating} on $K$ if for every $x,y\in K$ there is $U\in\Yu$ with $|U\cap\dn xy|=1$.

A compact space $K$ is {\em Valdivia compact} \cite{DG} if for some $\kappa$, $K$ embeds into $\Err^\kappa$ so that $\Sigma(\kappa)$ is dense in the image. Given $f\in\Err^\kappa$ and $S\subs\kappa$ denote by $f\mid S$ the function $(f\rest S)\concat 0_{\kappa\setminus S}$, i.e. $(f\mid S)(\al)=f(\al)$ if $\al\in S$ and $(f\mid S)(\al)=0$ otherwise. Let us recall an important factorization property of Valdivia compact spaces (see e.g. \cite{KM}): assuming $K\subs\Err^\kappa$ is such that $\Sigma(\kappa)\cap K$ is dense in $K$, for every map $\map fKY$ there exists $S\subs\kappa$ such that $|S|\loe\w(Y)$ and $f(x)=f(x\mid S)$ for every $x\in K$. In particular, $x\mid S\in K$ whenever $x\in K$. Applying this result to projections $\map{\pr_T}K{\Err^T}$ we see that $K$ is the limit of a continuous inverse sequence of smaller Valdivia compacta whose all bonding mappings are retractions.

\section{Semi-Eberlein compacta}

A compact space $K$ will be called {\em semi-Eberlein} if for some set $S$ there is an embedding $K\subs \Err^S$ such that $c_0(S)\cap K$ is dense in $K$. Clearly, every adequate compact space is semi-Eberlein and every semi-Eberlein compact space is Valdivia compact. Since there are examples of adequate Corson compacta which are not Eberlein \cite{Ta1}, not every semi-Eberlein Corson compact space is Eberlein.
Various examples of adequate Corson compacta which are not even Gul'ko compacta one could find in \cite{Le1}. The linearly ordered space $\omega_1+1$ is not semi-Eberlein (see below), so not every Valdivia compact space is semi-Eberlein. 

Below are two simple characterizations of semi-Eberlein compacta, one in terms of families of open $F_\sig$ sets and the other one in terms of spaces of continuous functions.
Both results are similar to known analogous characterizations of Eberlein compacta. 

\begin{proposition}\label{szop33} Let $K$ be a compact space. The following properties are equivalent.
\begin{enumerate}
	\item[(a)] $K$ is semi-Eberlein.
	\item[(b)] There exists a $T_0$ separating collection $\Yu$ consisting of open $F_\sig$ subsets of $K$ such that $\Yu=\bigcup_{\ntr}\Yu_n$ and the set $$\setof{p\in K}{(\forall\;\ntr)\;|\setof{U\in\Yu_n}{p\in U}|<\aleph_0}$$ is dense in $K$.
	\item[(c)] There exist a set $S=\bigcup_{\ntr}S_n$ and an embedding $\map hK{[0,1]^S}$ such that the set of all $p\in K$ with the property that $S_n\cap \suppt(h(p))$ is finite for every $\ntr$, is dense in $K$.
\end{enumerate}
\end{proposition}

\begin{pf} (a)$\implies$(b) We assume $K\subs \Err^S$ and $c_0(S)\cap K$ is dense in $K$. For each $r\in\Qyu\setminus\sn0$ and $s\in S$ define
$$U_{s,r}=\setof{x\in K}{x(s)>r}\quad\text{if }\; r>0$$
and
$$U_{s,r}=\setof{x\in K}{x(s)<r}\quad\text{if }\; r<0.$$
Then $U_{s,r}$ is an open $F_\sig$ set and $\Yu=\setof{U_{s,r}}{s\in S,\; r\in\Qyu\setminus\sn0}$ is $T_0$ separating on $K$. Finally, setting $\Yu_r=\setof{U_{s,r}}{s\in S}$, we have $\Yu=\bigcup_{r\in\Qyu\setminus\sn0}\Yu_r$ and for every $p\in c_0(S)\cap K$ there are only finitely many $s\in S$ with $|p(s)|>|r|$, i.e. the set $\setof{U\in\Yu_r}{p\in U}$ is finite.

(b)$\implies$(c) For each $U\in \Yu$ fix a continuous function $\map{f_U}K{[0,1]}$ such that $U=\inv f{(0,1]}$. Let $h$ be the diagonal product of $\sett{f_U}{U\in\Yu}$. Then $\map hK{[0,1]^\Yu}$ is such that $U\in\suppt(h(p))$ iff $p\in U$; therefore $\suppt(h(p))\cap \Yu_n$ is finite for every $\ntr$ and for every $p\in D$, where $D$ is the dense set of points described in (b).

(c)$\implies$(a) We assume $K\subs[0,1]^S$ and $S=\bigcup_{\ntr}S_n$ is as in (c). 
Let $D$ consist of all $p\in K$ such that $\suppt(p)\cap S_n$ is finite for every $\ntr$. Fix $s\in S$ and let $n$ be minimal such that $s\in S_n$. Define $\map{h_s}{[0,1]}{[0,1]}$ by $h_s(t)=t/n$. Then $h=\prod_{s\in S}h_s$ is an embedding of $[0,1]^S$ into itself. Now, if $p\in D$ and $\eps>0$ then $h(p)(s)>\eps$ only if $s\in S_n$ and $1/n>\eps$; thus $\setof{s\in S}{|h(p)(s)|>\eps}$ is finite. Hence $c_0(S)\cap \img hK$ is dense $\img hK$.
\end{pf}

We denote by $A(\kappa)$ the one-point compactification of the discrete space of size $\kappa$. Given a dense subset $D$ of a space $K$ we denote by $\tau_p(D)$ the topology of pointwise convergence on $D$ on the space $C(K)$.

\begin{proposition}\label{szop1} Let $K$ be a compact space. Then $K$ is semi-Eberlein iff there exist a dense set $D\subs K$ and a copy of $A(\kappa)$ in
 $\pair{C(K)}{\tau_p(D)}$ which separates the points of $K$. \end{proposition}

\begin{pf} Assume $K$ is semi-Eberlein, $K\subs \Err^\kappa$ and $D=K\cap c_0(\kappa)$ is dense in $K$. Let $$L=\setof{\pr_\al\rest K}{\al<\kappa}\cup\sn0,$$ where $\map{\pr_\al}{\Err^\kappa}{\Err}$ is the $\al$-th projection. We may assume that $\pr_\al\rest K\ne0$ for every $\al<\kappa$, shrinking $\kappa$ if necessary. We claim that $\pair L{\tau_p(D)}$ is homeomorphic to $A(\kappa)$ and $0$ is the only accumulation point of $L$. Fix a basic neighborhood $V$ of $0$ in $\pair{C(K)}{\tau_p(D)}$. Then
$$V=\setof{x\in C(K)}{(\forall\;i<k)\;|x(d_i)|<\eps}$$
for some $d_0,\dots,d_{k-1}\in D$ and $\eps>0$. 
By the fact that $d_i\in c_0(\kappa)$ for $i<k$, there is a finite set $S\subs \kappa$ such that $|d_i(\al)|<\eps$ for $i<k$ and for every $\al\in\kappa\setminus S$. Since $\pr_\al(d_i)=d_i(\al)$, we see that $\pr_\al\rest K\in V$ for $\al\in\kappa\setminus S$.

Now assume that $A(\kappa)$ embeds into $\pair{C(K)}{\tau_p(D)}$ for some dense set $D\subs K$ and its image separates the points of $K$. We may assume that $L=\setof{f_\al}{\al<\kappa}\subs C(K)$ is such that $L\cup\sn0$ endowed with the topology $\tau_p(D)$ is homeomorphic to $A(\kappa)$, $0$ is the accumulation point of $L$ and $L$ separates the points of $K$. 
Denote by $\En$ the (countable) collection of all nonempty open rational intervals $v\subs\Err$ such that $0\notin\cl v$. For each $v\in \En$ define
$U^v_\al=\inv{f_\al}v$. Let $\Yu_v=\setof{U^v_\al}{\al<\kappa}$. Then $\Yu=\bigcup_{v\in\En}\Yu_v$ is $T_0$ separating and consists of open $F_\sig$ sets. 
Fix $x\in D$ and $v\in \En$. There is a finite set $S\subs\kappa$ such that $f_\al(x)\notin\cl v$ for $\al\in\kappa\setminus S$ (because $\setof{f\in C(K)}{f(x)\in\Err\setminus\cl v}$ is a neighborhood of $0$ in $\tau_p(D)$).
If $\al\notin S$ then $x\notin U^v_\al$. Thus each $\Yu_v$ is point-finite on $D$. By Proposition \ref{szop33}, $K$ is semi-Eberlein.
\end{pf}

\begin{remark}
Every compact space $L$ can be embedded into
\newline 
$\pair{C(K)}{\tau_p(D)}$ for some compact space $K$ and a dense set $D\subs K$, in such a way that the image of $L$ separates the points of $K$. Indeed, assuming $L\subs[0,1]^\kappa$, consider first $K_0=\beta\kappa$, the \v Cech-Stone compactification of the discrete space $\kappa$, and let $D_0=\kappa$. Then $[0,1]^\kappa$ embeds naturally into $\pair{C(K_0)}{\tau_p(D_0)}$, since every bounded function on $\kappa$ extends
uniquely to $K_0$. Now define an equivalence relation $\sim$ on $K_0$ by
$$x\sim x' \iff (\forall\;f\in L)\;f(x)=f(x').$$ Let $K$ be the quotient space and let $\map q{K_0}{K}$ be the quotient map. Then $K$ is Hausdorff and $L$ naturally embeds into $\pair{C(K)}{\tau_p(D)}$, where $D=\img q{D_0}$. Clearly, after this embedding $L$ separates the points of $K$.
\end{remark}

\begin{remark}
It is known that if $K$ is a compact such that there exists a compact $L\subs C_p(K)$ which separates points of $K$ then $K$ is an 
Eberlein compactum.
This fact has no generalization to semi-Eberlein compacta: we give an example which shows that in Proposition \ref{szop1}
the compact $A(\kappa)$ cannot be replaced
by any other (even metrizable) compact space.
Take $K$ to be the \v Cech-Stone compactification of $\omega$ and let
$D=\omega$. Then the Cantor set $L=\{0,1\}^\omega$ is embedded naturally into
$\pair{C(K)}{\tau_p(D)}$. Clearly, this Cantor set $L$ separates the points of $K$.
$\beta\omega$ is not a continuous image of any Valdivia compact space \cite{Kalenda}.
\end{remark}

\section{A preservation theorem}

It is clear that the class of semi-Eberlein compacta is stable under arbitrary products and under closed subsets which have a dense interior. It is not hard to see that the one-point compactification of any topological sum of semi-Eberlein compacta is semi-Eberlein.
By \cite{Kalenda4} every non-Corson Valdivia compact space has a two-to-one map onto a non-Valdivia space. Thus, the class of semi-Eberlein compacta is not stable under continuous images.

It has been shown in \cite{KU} that there exists a compact connected Abelian topological group of weight $\aleph_1$ which is not Valdivia compact. As every compact group is an epimorphic (and therefore open) image of a product of compact metric groups, this shows that there exists a semi-Eberlein compact space (namely some product of compact metric groups) which has an open map onto a non-Valdivia space. On the other hand, it has been shown in \cite{KM} that ``small" 0-dimensional open images as well as retracts of Valdivia compact spaces are Valdivia, where ``small" means ``of weight $\loe\aleph_1$". Another result from \cite{KM} states that the class of Valdivia compacta is stable under limits of certain inverse sequences with retractions \cite{KM}; in particular a limit of a continuous inverse sequence of metric compact spaces whose all bonding mappings are retractions is Valdivia compact (in fact this characterizes Valdivia compact spaces of weight $\loe\aleph_1$). 
Below we prove a parallel preservation theorem for semi-Eberlein spaces. As an application we show that retracts of Cantor or Tikhonov cubes are semi-Eberlein, which improves a result from \cite{KM}. Note that $\omega_1+1$ is an example of a non-semi-Eberlein space which is the limit of a continuous inverse sequence of metric compacta whose all bonding mappings are retractions.

Recall that a map $\map fXY$ is called {\em semi-open} if for every nonempty open set $U\subs X$ the image $\img fU$ has nonempty interior.

\begin{lemma}\label{gora}
Assume $\map fXY$ is a retraction, i.e. $Y\subs X$ and $f\rest Y=\id_Y$. Assume $\map gX{[0,1]^\omega}$ is such that $f\diag g$ is one-to-one. Then $Y$ is $G_\delta$ in $X$ and for every $T_0$ separating family $\Yu$ on $Y$ there exists a countable family $\Vee$ of open $F_\sig$ sets in $X$ such that $\Yu\cup\Vee$ is $T_0$ separating on $X$ and $Y\cap \bigcup\Vee=\emptyset$.
\end{lemma}

\begin{pf}
Fix a metric $d$ on $[0,1]^\omega$ and define $\phi(x)=d(g(x),g(f(x)))$. Then $\phi(x)=0$ iff $g(x)=g(f(x))$ which is equivalent to $x=f(x)$, by the fact that $f\diag g$ is one-to-one. Thus $Y=\phi^{-1}(0)$ is $G_\delta$ in $X$. Now, fix a countable $T_0$ separating open family $\Vee_0$ on $[0,1]^\omega$ and define $\Vee=\setof{\inv gV\setminus Y}{V\in \Vee_0}$. Clearly, $\Vee$ consists of open $F_\sig$ sets and if $g(x)\ne g(x')$ then there is $V\in \Vee$ such that $|\{x,x'\}\cap V|=1$. Since $f\diag g$ is one-to-one, $\Yu\cup \Vee$ is $T_0$ separating on $X$.
\end{pf}

A map $\map fXY$ has a {\em metrizable kernel} if there exists a map $\map gX{[0,1]^\omega}$ such that $f\diag g$ is one-to-one.

\begin{theorem}\label{krater} Assume $\S=\invsys Kr\kappa \al\beta$ is a continuous inverse sequence such that
\begin{enumerate}
	\item $K_0$ is semi-Eberlein,
	\item each $r^{\al+1}_\al$ is a semi-open retraction with a metrizable kernel.
\end{enumerate}
Then $\liminv\S$ is semi-Eberlein.
\end{theorem}

\begin{pf}
Let $K=\liminv\S$. We may assume that each $K_\al$ is a subspace of $K$ and each projection $\map{r_\al}K{K_\al}$ is a retraction (see \cite[Prop. 4.6]{BKT}); therefore $r^\beta_\al=r_\al\rest K_\beta$ for $\al<\beta$. We construct inductively families of open $F_\sig$ sets $\Yu_\al=\bigcup_{\ntr}\Yu_\al^n$ and dense sets $D_\al\subs K_\al$ such that:
\begin{enumerate}
  \item[(a)] $\Yu_\al^n$ is point-finite on $D_\al\subs K_\al$;
  \item[(a')] each $U\in \Yu_\al$ is of the form $\inv{r_\al}{U'}$ for some open $F_\sig$ subset of $K_\al$;
	\item[(b)] $\al<\beta\implies \Yu_\al^n\subs\Yu_\beta^n$ and $K_\al\cap\bigcup(\Yu_\beta^n\setminus\Yu_\al^n)=\emptyset$ for every $\ntr$;
	\item[(c)] if $\gamma$ is a limit ordinal then $\Yu_\gamma^n=\bigcup_{\al<\gamma}\Yu_\al^n$ and $D_\gamma=\bigcup_{\al<\gamma}D_\al$;
	\item[(d)] $D_{\al+1}=\inv{(r^{\al+1}_\al)}{D_\al}$.
	\item[(e)] $\Yu_\al$ is $T_0$ separating on $K_\al$. 
\end{enumerate}
We start by finding a suitable family $\Vee_0=\bigcup_{\ntr}\Vee_0^n$ on $K_0$, using Proposition \ref{szop33}(b). Let $D_0$ be a dense set with the property that every $\Vee_0^n$ is point-finite on $D_0$. Define $\Yu_0^n=\sett{\inv{r_0}V}{V\in\Vee_0^n}$. Clearly (a), (a') and (e) hold.

The successor stage is taken care of by Lemma \ref{gora}: we use the fact that $r^{\al+1}_\al$ is semi-open to deduce that $\inv{(r^{\al+1}_\al)}{D_\al}$ is dense in $K_{\al+1}$. Note also that if $V\subs K_{\al+1}$ is disjoint from $K_\al$ then $\inv{r_{\al+1}}V\cap K_\al=\emptyset$, which ensures that condition (b) is satisfied.

It remains to check the limit stage. Fix a limit ordinal $\delta\loe\kappa$. Clearly, $D_\delta$ defined by (c) is dense and $\Yu_\delta$ defined in (c) consists of open $F_\sig$ sets and satisfies (a'). Fix $x\in D_\delta$ and $\ntr$. Then $x\in D_\al$ for some $\al<\delta$ and thus by induction hypothesis $x$ belongs to only finitely many elements of $\Yu_\al^n$. By the second part of (b), $x\notin U$ for any $U\in\Yu_\beta^n\setminus \Yu_\al^n$. It follows that each $\Yu_\delta^n$ is point-finite on $D_\delta$. 

Now fix $x\ne x'$ in $K_\delta$. Then $r_\al(x)\ne r_\al(x')$ for some $\al<\delta$ (by the continuity of the sequence) and hence there is $U\in\Yu_\al$ such that e.g. $r_\al(x)\in U$ and $r_\al(x')\notin U$. By (a') we have $U=\inv{r_\al}{U'}$. Note that $U'=U\cap K_\al$, because $r_\al\rest K_\al=\id_{K_\al}$. Thus $r_\al(x)\in U'$, $r_\al(x')\notin U'$ and hence $x\in U$ and $x'\notin U$. It follows that $\Yu_\delta$ is $T_0$ separating. 

Thus the construction can be carried out. By Proposition \ref{szop33}, each $K_\al$, and therefore also $K$, is semi-Eberlein. 
\end{pf}

\begin{corollary}
Every retract of a Cantor or Tikhonov cube is semi-Eberlein.
\end{corollary}

\begin{pf}
Let $K$ be a retract of a Cantor or Tikhonov cube. Then $K=\liminv\S$, where $\S=\invsys Kr\kappa\al\beta$ is a continuous inverse sequence such that $K_0$ is a compact metric space and each $r^{\al+1}_\al$ is an open retraction with a metrizable kernel. In both cases Theorem \ref{krater} applies.

In the case of Cantor cubes this is Haydon--Koppelberg decomposition theorem \cite[Corollary 2.8]{Kop}. In fact this result is about {\em projective Boolean algebras} which, via Stone duality, correspond to retracts of Cantor cubes.

In the case of Tikhonov cubes, one needs to refer to Chapter 2 of Schepin's work \cite{Szcz}. Since the quoted result is not proved explicitly in \cite{Szcz} and we could not find a suitable bibliographic reference for it, we briefly sketch the proof, giving references to appropriate results from \cite{Szcz}. For a full and self-contained proof we refer to \cite{K}.

Assume $\map r{[0,1]^\kappa}K$ is a retraction. Following Shchepin, we say that a set $S\subs \kappa$ is {\em $r$-admissible} if $x\rest S=x'\rest S$ implies $r(x)\rest S=r(x')\rest S$ for every $x,x'\in[0,1]^\kappa$. This property is equivalent to the following: if $f\in C(K)$ depends on $S$, i.e. $x\rest S=x'\rest S\implies f(x)=f(x')$, then $fr\in C([0,1]^\kappa)$ also depends on $S$. Define $K_S=\img{\pr_S}K$ and let $\map {p_S}K{K_S}$ be defined by $p_S=\pr_S\rest K$, where $\map{\pr_S}{[0,1]^\kappa}{[0,1]^S}$ is the projection. For every $r$-admissible set $S\subs\kappa$ the map $p_S$ is open (Lemma 1 in \cite[Chapter 2]{Szcz}) and soft (Lemma 6 in \cite[Chapter 2]{Szcz}), thus in particular it is a retraction. Observe that the union of any family of $r$-admissible sets is $r$-admissible. By Lemma 3 in \cite[Chapter 2]{Szcz}, every countable subset of $\kappa$ can be enlarged to a countable $r$-admissible set. For each $\al<\kappa$ choose a countable $r$-admissible set $S_\al\subs\kappa$ such that $\al\in S_\al$. Define $A_\al=\bigcup_{\xi<\al}S_\xi$. Then each $A_\al$ is $r$-admissible and hence each $p_\al:=p_{A_\al}$ is an open soft map. Given $\al\loe \beta<\kappa$, let $p^\beta_\al$ be the unique map such that $p_{\al}= p^\beta_\al p_\beta$. Then $\Es=\invsys Kp\kappa\al\beta$, where $K_\al=K_{A_\al}$, is a continuous inverse sequence such that each $p^{\al+1}\al$ is an open retraction with a metrizable kernel (since $S_\al$ is countable) and $K=\liminv \Es$.
\end{pf}

\section{Semi-Eberlein spaces have no P-points}

Recall that a point $p\in X$ is called a {\em P-point} if $p$ is not isolated and $p\in\Int\bigcap_{\ntr}U_n$ for every sequence $\setof{U_n}{\ntr}$ of neighborhoods of $p$. The simplest example of a compact space with a P-point is $\omega_1+1$, considered as a linearly ordered space. In this section we observe that semi-Eberlein spaces do not have P-points, which implies in particular that $\omega_1+1$ is not semi-Eberlein. Using this observation and a simple forcing argument, we show that the Corson compact of Todor\v cevi\'c defined in \cite{To2} is not semi-Eberlein (Example \ref{bajer} below).

A net $\sett{x_\sig}{\Ssig}$ in a topological space $X$ will be called {\em $\sig$-bounded} if $\Sigma$ is $\sig$-directed and for every countable set $S\subs\Sigma$ there is $\tau\in\Sigma$ such that $\cl\setof{x_\sig}{\sig\goe \tau}\cap \cl\setof{x_\sig}{\sig\in S}=\emptyset$. In case where $\sett{x_\sig}{\Ssig}$ converges to $p\in X$, this just means that $p\notin\cl\setof{x_\sig}{\sig\in S}$ for any countable $S\subs\Sigma$.

\begin{lemma}\label{sarna1} For every set $S$, there are no $\sig$-bounded nets in $c_0(S)$ which converge in $\Err^S$.
\end{lemma}

\begin{pf} Suppose $\sett{x_\sig}{\Ssig}\subs c_0(S)$ is a $\sig$-bounded net which converges to $p\in\Err^S$. Then $\suppt(p)$ is uncountable, because
\newline
 $\setof{x_\sig}{\Ssig}\subs\Sigma(S)$ and $\Sigma$-products of real lines have countable tightness.
Choose $\setof{s_n}{\ntr}\subs\suppt(p)$ so that $|p(s_n)|>\eps$ for some fixed $\eps>0$. Using the fact that $p=\lim_{\Ssig}x_\sig$, for each $\ntr$ we find $\sig_n\in\Sigma$ such that $|x_\sig(s_n)|>\eps$ for every $\sig\goe\sig_n$. Let $\tau\in\Sigma$ be such that $\sig_n<\tau$ for every $\ntr$. Then $|x_\tau(s_n)|>\eps$ for every $\ntr$ and consequently $x_\tau\notin c_0(S)$, a contradiction.
\end{pf}

\begin{theorem}\label{hop1} Semi-Eberlein compact spaces do not have P-points. \end{theorem}

\begin{pf} Assume $K\subs\Err^\kappa$, $c_0(\kappa)\cap K$ is dense in $K$ and suppose that $p\in K$ is a P-point. Then $p\notin c_0(\kappa)$. Let $\Sigma$ be a fixed base at $p$. Consider $\Sigma$ as a ($\sig$-directed) poset with reversed inclusion. For each $u\in\Sigma$ choose $x_u\in u\cap c_0(\kappa)$. If $S\subs\Sigma$ is countable then there is $v\in\Sigma$ such that $v\subs K\setminus\setof{x_u}{u\in S}$ and therefore $p\notin \cl\setof{x_u}{u\in S}$. Thus $\sett{x_u}{u\in\Sigma}$ is a $\sig$-bounded net in $c_0(\kappa)$ which converges to $p$. By Lemma \ref{sarna1}, we get a contradiction.
\end{pf}

\begin{corollary} The linearly ordered space $\omega_1+1$ is not semi-Eberlein. \end{corollary}

Actually, we can strengthen Theorem \ref{hop1} by saying that if $K$ is semi-Eberlein then no forcing notion can force a P-point in $K$. More precisely, if $K$ is a compact space, $\poset$ is a forcing notion and $G$ is a $\poset$-generic filter then by $K^G$ we denote the compact space in the $G$-generic extension which consists of all ultrafilters in $CL(K)$, where $CL(K)$ denotes the lattice of closed subsets of $K$ defined in the ground model. It is not hard to check that we get the same space if we take, instead of $CL(K)$, any of its sublattices which is a closed base for $K$. Indeed, if $p\ne q$ are ultrafilters in $CL(K)$ and $\bbl$ is a sublattice of $CL(K)$ which is a closed base then, by compactness, there are disjoint $a,b\in\bbl$ such that $a\in p$ and $b\in q$. This shows that the map $\map h{\ult(CL(K))}{\ult(\bbl)}$, defined by $h(p)=p\cap \bbl$, is one-to-one ($\ult(\bbl)$ denotes the space of all ultrafilters over $\bbl$). Clearly, $h$ is a continuous surjection and, when defined in a fixed generic extension, it shows that $\ult(\bbl)$ is homeomorphic to $\ult(CL(K))$. 

It follows that $[0,1]^G$ is the same as the usual unit interval in the $G$-generic extension, since it can be described as the space of ultrafilters over the lattice generated by all closed rational intervals. A similar fact is true for any cube $[0,1]^\kappa$.

Next, observe that $K$ is dense in $K^G$. Indeed, $K\subs K^G$, because being an ultrafilter in a lattice is absolute. A basic open set in $K^G$ is of the form $a^-=\setof{p\in K^G}{a\notin p}$, where $a$ is an element of $CL(K)$ from the ground model. Now, if $a^-\nnempty$ then there is $b\in CL(K)$ such that $a\cap b=\emptyset$ and therefore in the ground model there is $p\in \ult(CL(K))$ such that $b\in p$; consequently $K\cap a^-\nnempty$.

Finally, if $K,L$ are compact spaces such that $K\subs L$ then $K^G$ is the closure of $K$ in $L^G$. Indeed, the fact $K\subs L$ is encoded in the lattice epimorphism $\map h{\bbl}{\bbk}$, defined by $h(a)=a\cap K$, where $\bbk=CL(K)$ and $\bbl=CL(L)$. In any generic extension, $h$ induces an embedding of $K^G$ into $L^G$ and since $K$ is dense in $K^G$ and $K^G$ is closed, we have $K^G=\cl_{L^G}(K)$.

In particular, if $K\subs[0,1]^\kappa$ then $K^G$ can be regarded as the closure of $K$ in $([0,1]^\kappa)^G$, where $([0,1]^\kappa)^G$ is the Tikhonov cube $[0,1]^\kappa$ defined in the $G$-extension. Clearly, $K^G$ is semi-Eberlein if $K$ is so. Thus:

\begin{proposition} No forcing notion can introduce a P-point in a semi-Eberlein space. More precisely, if $K$ is a semi-Eberlein space, $\poset$ is a forcing notion and $G$ is a $\poset$-generic filter then $K^G$ has no P-points. \end{proposition}

This observation leads to examples of Corson compact spaces which are not semi-Eberlein.

\begin{example}[Todor\v cevi\'c]\label{bajer}
There exists a Corson compact space
which is not semi-Eberlein. In fact the Corson
compact space of \cite{To2} (see also \cite[p. 287]{To1}) is such an example.
Let us recall the construction of this space. For a tree $T$,
let $P(T)$ be the set of all initial branches of $T$ (i.e. linearly ordered sets $x\subs T$ such that $t\in x\Land s<t\implies s\in x$) with the
topology induced from the Cantor cube $\{0,1\}^T$. When $T$ has no uncountable 
branches, $P(T)$ is an example of a Corson compact space. The tree $T$
used in \cite{To2} is the set of all subsets of
some fixed stationary and costationary set $A\subs \omega_1$ which are closed in $\omega_1$.
Thus, in particular, such a tree $T$ is Baire (see \cite[Lemma 9.12]{To1})
and therefore forcing with $T$ does not collapse $\omega_1$ (see \cite[Theorem 15.6]{Je}, where a Baire forcing notion is called $\aleph_0$-distributive).
A $T$-generic filter $G\subs T$ gives an uncountable branch which is obviously a P-point in $P(T)^G$. Hence, $P(T)$ is not semi-Eberlein.
\end{example}

\begin{remark}
It is shown in \cite{To2} that the Corson compactum $P(T)$ has no dense metrizable subspaces.
Evidently, for each tree the family of all its chains is an adequate family.
Note that the adequate (and hence semi-Eberlein) Corson compactum built on this tree $T$
also has no dense metrizable subspaces \cite{Le2}. This shows that a semi-Eberlein compact which is Corson does not have to be Eberlein or even Gul'ko compact (since Gul'ko compacta contain dense completely metrizable subspaces \cite{G}). 
\end{remark}

\section{Questions}

There are several open questions concerning Valdivia compacta; one may ask parallel questions for semi-Eberlein spaces. Below are questions specific for semi-Eberlein compacta.

\begin{question}
Can a semi-Eberlein compact space have weak P-points?
\end{question}

\begin{question}
Is the class of semi-Eberlein spaces stable under closed $G_\delta$ sets?
\end{question}

Even for adequate compacts the same question seems to be open.
\begin{question}
Let $X$ be an adequate compact and $F\subs X$ is a closed $G_\delta$ set.
Is $F$ a semi-Eberlein compact space?
\end{question}
As a relevant result, we mention that  Bell \cite{Be} gave an example of a centered (continuous image of adequate) compact $X$ and 
a closed $G_\delta$ set $Z\subs X$ 
such that $Z$ is not homeomorphic to any centered compact. However $Z$ in his example is an Eberlein compact. 

\begin{question} Does there exist a semi-Eberlein space $K$ which does not embed into any cube $[0,1]^\kappa$ so that $\sig(\kappa)\cap K$ is dense in $K$? ($\sig(\kappa)$ is the small sigma-product of $\kappa$ copies of $\Err$)
\end{question}

\begin{question}
Assume $K$ is semi-Eberlein, not Eberlein. Does $K$ have a non-semi-Eberlein image?
\end{question}

\begin{question}
Assume $K$ is semi-Eberlein and $\map fKL$ is a retraction or an open surjection such that $L$ has densely many $G_\delta$ points. Is $L$ semi-Eberlein?
\end{question}

\bibliographystyle{amsplain}

\end{document}